\documentclass[11pt]{article}
\usepackage{amsmath}
\usepackage{amssymb}
\usepackage{url}
\usepackage{booktabs}
\newcommand{\setcard}[1]{\ensuremath{\left| #1 \right|}}

\author{Drago\v s Cvetkovi\' c, Jason Grout}

\title{Graphs with extremal energy should have a small number of
  distinct eigenvalues}
\begin{document}
\maketitle

\begin{abstract}

  {\em The sum of the absolute values of the eigenvalues of a graph is
    called the energy of the graph. We study the problem of finding
    graphs with extremal energy within specified classes of graphs. We
    develop tools for treating such problems and obtain some partial
    results.  Using calculus, we show that an extremal graph
    ``should'' have a small number of distinct eigenvalues.  However,
    we also present data that shows in many cases that extremal
    graphs can have a large number of distinct
    eigenvalues.\footnote{The idea for this note arose during the
      workshop ``Spectra of families of matrices described by graphs,
      digraphs, and sign patterns'' which was held at American
      Institute of Mathematics in Palo Alto, California, U.S.A. on
      October 23--27, 2006.}}

\end{abstract}



\section{Introduction}

Let $G$ be a simple undirected graph on $n$ vertices with $m>0$ edges. Let
$x_1,x_2,\ldots,x_n$ be the eigenvalues of $G$. Let
$I=\{1,2,\ldots,n\}$.

The following relations are well-known:

\begin{align}
\sum_{i\in I}x_i=0,\label{eq:first-moment}\\
\sum_{i\in I}x_i^2=2m\label{eq:second-moment}.
\end{align}

The energy $E(G)$ of $G$ is defined by
$$E(G)=\sum_{i\in I}|x_i|.$$
The energy of a graph was defined by I. Gutman \cite{Gut1} and has
attracted much attention from researchers in the last few years.

Given a set $\cal G$ of graphs, one can ask which graphs in $\cal G$
have extremal (minimal or maximal) energy.  In this paper, we shall
develop a procedure for treating the problem of finding an extremal
graph when all graphs in $\cal G$ have a fixed number of vertices and
edges and a given family of eigenvalues.  Because our procedure allows
us to specify eigenvalues that graphs in $\cal G$ must have, we can
make sure that $\cal G$ includes all graphs having certain properties
(bipartite, regular, etc.).  Because our procedure
relies on methods of continuous calculus, it will not always produce a
graph with extremal energy.  However, even when we do not obtain a
graph with extremal energy, our procedure still gives a heuristic to
guide a search for such a graph.

We will concentrate our efforts in this paper on finding a graph in
$\cal G$ that has maximal energy.  Our tools can also be applied to
the problem of finding graphs with a minimal energy. However, this
problem appears to be easier and there are several recent results in
this direction (see, for example, \cite{Yan}, \cite{Hua}), so we shall
mention it only in passing.

It is known that for $n\leq 7$, the graphs with maximal energy are the
complete graphs $K_n$, $n=1,2,\dots,7$. Maximal values of energy for
graphs with $n$ vertices have been determined heuristically by the system
AutoGraphiX \cite{CCGH} for $n\leq 12$.  In
\cite{CCGH}, the corresponding maximal energy graph is given only for
$n=10$.  The graph with maximal
energy among 10-vertex graphs is the complement of the Petersen graph,
i.e., the line graph $L(K_5)$ of the complete graph $K_5$.  This graph
is strongly regular and has three distinct eigenvalues. Maximal energy
graphs have been determined in \cite{KM} for an infinite sequence of
values of $n$ ($n$ is a power of
2) and in these cases the graphs are also strongly regular. The
smallest such graph, the Clebsch graph, has 16 vertices and the
spectrum $10, 2^5, -2^{10}$ (superscripts denote multiplicities of
eigenvalues). The Clebsch graph appears as graph
No.~187 in Table~A3 in~\cite{CvRS4}.

The connection between maximal energy graphs and strong regularity is
explored further in \cite{KM}.  It is proved in \cite{KM} that for a
graph $G$ on $n$ vertices the following inequality holds:
$$ E(G) \leq \frac {n}{2} (1+ \sqrt{n}), $$
with equality if and only if $G$ is a strongly regular graph with
parameters $(n,(n+\sqrt{n})/2,(n+\sqrt{n})/4,(n+\sqrt{n})/4)$.
Such strongly regular graphs exist for $n=4\tau ^2$ for  $\tau =
2^m$, $m=1,2,\dots$. From the theory of strongly regular graphs, one
can deduce that a graph with such parameters has distinct
eigenvalues $\tau(2\tau +1), \pm \tau$. The conjecture that, for
any $\epsilon >0$, for almost all $n$, there exists a graph $G$ on $n$
vertices such that
$$ E(G) \geq (1-\epsilon)\frac {n}{2} (1+ \sqrt{n}) $$
has been confirmed in a slightly improved form in
\cite{Nik}.
These results give a basis to think that, in general, graphs with
extremal energy have a small number of distinct eigenvalues. Though
our results support this heuristic, we also note examples with
significant deviations from these expectations.  These unexpected
examples seem to make the maximal energy problem very difficult.

The rest of the paper is organized as follows. In Section~2, we
present some known results on graph spectra which will be used later.
In Section~3, we prove a simple but useful theorem which gives a basis
for thinking that extremal graphs should have a small number of
distinct eigenvalues.  Section~4 contains examples using our theorem.
In Section~5, we present some computational and theoretical data which
show the existence of maximal energy graphs with a large number of
distinct eigenvalues. A heuristic procedure for finding maximal energy
graphs is formulated in Section~6. Finally, we summarize our main
points concerning maximal energy graphs in Section~7.

\section{Preliminaries}

We shall quote some simple and well-known facts from the theory of
graph spectra which will be used in the rest of the paper.

\medskip\noindent
{\bf Proposition 1.} (see \cite{CvDSa}, p.~85) {\em  For any graph
with eigenvalues $x_1,x_2,\dots,x_n$ and $t$ triangles we have
$\sum_{i\in I}x_i^3=6t$.}

\medskip
In particular, this means that the third spectral moment
$\sum_{i\in I}x_i^3$ is always an integer that is divisible by 6.

\medskip\noindent
{\bf Proposition 2.} (see \cite{CvDSa}, p.~87) {\em A graph is
bipartite if and only if its spectrum is symmetric with respect to
0.}

\medskip
In bipartite graphs, if $x$ is an eigenvalue of multiplicity $k$,
then $-x$ is also an eigenvalue of multiplicity $k$.

\medskip\noindent
{\bf Proposition 3.} (see \cite{CvDSa}, p.~94) {\em A graph with
$n$ vertices, $m$ edges, and largest eigenvalue $x_1$ is
regular if and only if $x_1=\frac{2m}{n}$. }

\medskip\noindent {\bf Proposition 4.} (see \cite{CvDSa}, p.~56) {\em
  If $G$ is a regular graph with eigenvalues $x_1$,
  $x_2$, $\dots$, $x_n$ (with $x_1\geq x_i$ for $i=2,\ldots,n$), then the
  eigenvalues of the complement of $G$ are $n-x_1-1$ and $-x_i-1$ for
  $i=2,\ldots,n$.}

\medskip\noindent
{\bf Proposition 5.} (see \cite{CvDSa}, p.~53) {\em The eigenvalues of
a cycle $C_n$ of length $n$ are $2 \cos \left(\frac{2\pi}{n}
j\right)$, where $j=1,2,\dots ,n$. }

\medskip
In particular, we obtain the spectrum for a triangle ($2,-1^2$), a
quadrangle ($2,0^2,-2$), and a pentagon ($2, 0.618^2, -1.618^2$, i.e., 
$2, (1/\varphi)^2, -\varphi^2$,
where $\varphi=(1+\sqrt{5})/2 \approx 1.618$, the golden ratio).

\section{Main theorem and tools}
To motivate our main theorem, we shall use traditional calculus to try
to find extremal values of the energy of graphs. Of course, using
techniques from continuous mathematics in problems of discrete
mathematics requires careful handling and we shall see the limitations of
such an approach. However, some results can be achieved. The use of
calculus in handling extremal problems with graph eigenvalues was
suggested in \cite[Section~7.7]{CvDSa}.

Define
$$I=I_+\cup I_-, \text{ where }$$
$$I_+=\{i|i\in I,x_i\geq0\}\quad \text{and} \quad I_-=\{i|i\in I,x_i<0\}.$$
Then the energy can be represented in the following form:
$$E=\sum_{i\in I_+}x_i-\sum_{i\in I_-}x_i.$$
If the graph is not trivial, then $I_+$ and $I_-$ are
both non-empty by (\ref{eq:first-moment}).

Consider an auxiliary function involving the constraints
(\ref{eq:first-moment}) and (\ref{eq:second-moment}):
$$F=\sum_{i\in I_+}x_i-\sum_{i\in I_-}x_i +
\alpha\sum_{i\in I}x_i+\beta\left(\sum_{i\in
I}x_i^2-2m\right),$$
where $\alpha,\beta$ are Lagrange multipliers.
Extremal values of the function $E$ satisfying the constraints
(\ref{eq:first-moment}) and (\ref{eq:second-moment}) can be found by
equating the partial derivatives of function $F$ with 0, i.e.
$$\frac{\partial F}{\partial
x_j}=\pm 1+\alpha+2\beta x_j=0,~~j \in I.$$
The first term in the sum is equal to $+1$ if $j \in I_+$ and is equal
to $-1$ if $j \in I_-$. Hence we obtain
$$x_j=\frac{-\alpha\mp
  1}{2\beta},~~j \in I.$$ This means that a graph with extremal energy
should have only two distinct eigenvalues. 

These derivations were used in \cite{Gut} to derive the McClelland
upper bound for the energy in a new way.

However, the only graphs that have only two distinct eigenvalues are
unions of complete graphs of a fixed order.  We must assume that at
least three distinct eigenvalues exist in non-trivial cases.

We shall treat our problem in the following way. Assume that we
somehow know that an extremal graph we are looking for has some fixed
and given eigenvalues $\cal K$, consisting of $x_i$, $i\in K$. Let
${\cal H}$ be the set of graphs having a fixed number of vertices and
edges and the given family of eigenvalues $\cal{K}$.  Let
$J=I\setminus K$, so that the eigenvalues $x_i$, $i\in J$, are
considered unknown.  We seek to determine these unknown eigenvalues in
such a way that the energy becomes extremal in $\cal{H}$.  Index sets
$J$ and $K$ are further partitioned, as $I$ was partitioned, into
subsets corresponding to non-negative and negative eigenvalues:
$$J=J_+\cup J_-
\quad \text{and} \quad K=K_+\cup K_-.$$

Now we have
$$E=\sum_{i\in J_+}x_i-\sum_{i\in
J_-}x_i+\sum_{i\in K_+}x_i-\sum_{i\in K_-}x_i,$$
and from (\ref{eq:first-moment}) and (\ref{eq:second-moment}),
$$\sum_{i\in
J}x_i+\sum_{i\in K}x_i=0\quad \text{and} \quad \sum_{i\in J}x_i^2+\sum_{i\in
K}x_i^2=2m.$$
Let 
\begin{align}
C_+=\sum_{i\in K_+}x_i, \; C_-=\sum_{i\in K_-}x_i, \; C=\sum_{i\in K}x_i, \text{ and } D=\sum_{i\in K}x_i^2.\label{eq:CD}
\end{align}
We can write
$$F=\sum_{i\in J_+}x_i-\sum_{i\in J_-}x_i
+C_+-C_-+\alpha\left(\sum_{i\in
J}x_i+C\right)+\beta\left(\sum_{i\in J}x_i^2+D-2m\right).$$
Using partial derivatives for any $j \in J$, we get
$$\frac{\partial F}{\partial x_j}=\pm
1+\alpha+2\beta x_j=0 \implies x_j=\frac{-\alpha\mp 1}{2\beta}.$$
Assuming that both sets $J_+$ and $J_-$ are non-empty, this means that
unknown eigenvalues should have only two values in a graph with
extremal energy. (Sometimes one of the sets $J_+$ or $J_-$ is
 empty in which case our approach does not give any solution.)  Denote 
the two values obtained by $x,y$ and the corresponding multiplicities by $p,q$.

This sets up $|J|-1$ Lagrange multiplier problems (one for each situation
$|J_+|=i$ and $|J_-|=|J|-i$ for $i=1,2,\dots,|J|-1$).  For a given
distribution of unknown positive and negative eigenvalues, the
solution of the corresponding Lagrange multiplier problem will give us
an upper bound on the maximal energy of graphs in $\cal H$ with a
corresponding distribution of unknown eigenvalues.  If we take the
maximal value $E$ over all such solutions, and that energy value is
realized by a graph, then we know we have a maximal energy graph in
$\cal{H}$.

Extend the known part of the spectrum ${\cal K}$ by the values corresponding to $E$ of $x,y$ with multiplicities $p,q$. Denote by $L$ the described procedure of completing the partial spectrum using the Lagrange multipliers and let ${\cal K}_L$ be the resulting complete spectrum.

Now we can formulate a slightly more general statement.

\medskip
\noindent{\bf Theorem 1.} {\em Let $n,m$ be positive integers and let
  $\cal K$ be a family of reals.  Let $\cal G$ be the class of graphs
  having $n$ vertices and $m$ edges.   Suppose that $G\in{\cal G}$ is a graph with
maximal energy in $\cal G$ and  has all members of $\cal
  K$ as eigenvalues. Let ${\cal K}_L$ be the spectrum obtained
from $\cal K$ by procedure $L$. If a graph $H$ has the spectrum ${\cal K}_L$, 
then $H$ also has maximal energy in $\cal G$ and has the same
energy as $G$.}  \medskip

\noindent{\bf Proof.} Suppose that $G\in{\cal G}$ is a maximal energy graph in
$\cal G$ with energy $E$. Let ${\cal H}$ be the class of graphs from ${\cal G}$
which have all members of ${\cal K}$ as eigenvalues. By the construction of
${\cal K}_L$, the graph $H$ belongs to $\cal H$ and has a maximal energy in $\cal H$.
Since $G\in{\cal H}$, $H$ must have
energy at least $E$. However, since $G$ is also a maximal graph in ${\cal G}$,
the energy of $H$ must be exactly equal to $E$, so $H$ also is a maximal
graph in ${\cal G}$.\hfill \ensuremath{\Box}

\medskip

An analagous procedure and theorem can be stated for minimal energy graphs.

This theorem can be interpreted as ``graphs with extremal energy {\em
  should} have a small number of distinct eigenvalues,'' as pointed out
in the title of this paper. Namely, the eigenvalues which are not
fixed by the set $\cal K$ would have, in a truly optimal graph, only two
values, if such a graph happened to exist.

\section{Examples}

In this section we shall give several examples of using Theorem~1 to
guide searches for maximal energy graphs.  Note that these examples
may not represent complete proofs that the graphs involved have extremal
energy. However, they do illustrate how to use the ideas behind
Theorem~1.

As a computational alternative, instead of explictly using procedure
$L$ and the involved Lagrange equations, we will compute all possible
ways of adding two eigenvalues to the family $\cal K$ while respecting
the constraints (\ref{eq:first-moment}) and (\ref{eq:second-moment}).
The solutions of the Lagrange equations in procedure $L$ will be
among these results, so maximizing the energy over these results will
yield the maximum we would get from procedure $L$.

To this end, we introduce some technical notation.  The set $J$
contains $|J|=n-\setcard{K}$ elements. Suppose that an extremal solution
contains eigenvalues $x,y$ with multiplicities $p,q$ respectively.
The following system of equations must be satisfied, where $p$, $q$,
and $|J|$ are positive integers and $x$, $y$, $C$, and $D$ are real
numbers:
\begin{align}\{   p + q = |J|,\,
     px + qy = -C,\,
     px^2 + qy^2 = 2m-D
  \}.\label{eq:constraints}\end{align}
In this notation, the energy is $E = p|x| + q|y| + C_+ - C_-$ .

Our examples are computed using a small Mathematica program (see
Appendix~\ref{app:mathematica-program}) or a small SAGE \cite{sage} program (see
Appendix~\ref{app:sage-program}).  Given the number of vertices, the
number of edges, and a list of known eigenvalues that define $\cal G$,
the program prints all 4-tuples $(p,q,x,y)$ satisfying the system of
equations (\ref{eq:constraints}).  For each solution, the program also
prints the corresponding energy $E$ and $\frac{1}{6}\sum_{i\in I}x_i^3$.
This last quantity is the third spectral moment divided by 6 and
should be an integer if the solution corresponds to the eigenvalues of
an actual graph (see Proposition~1).  To solve the system of
equations, we vary $p$ between 1 and $|J|-1$ and for each such $p$
corresponding solutions are obtained (0, 1 or 2 solutions, having in
mind that the third equation in (\ref{eq:constraints}) is quadratic).
However, because of the symmetry of the system of equations, it is
sufficient to consider solutions for $p=1,2,\dots,\lfloor |J|/2 \rfloor$.

\medskip
\noindent {\bf A.}~ Let us examine the class of regular graphs of
degree 10 having 16 vertices. By Proposition~3, the class $\cal G$ defined by
$n=16$, $m=80$, and ${\cal K}=\{10\}$ coincides with this class.  This implies
that $\setcard{K}=1$, $|J|=15$, $C=10$, and $D = 100$.  Our program gives the
following results.

\begin{center}
\begin{tabular}{rrr@{.}lr@{.}lr@{.}lr@{.}l}\toprule
 \multicolumn{1}{c}{$p$} 
&\multicolumn{1}{c}{$q$}
&\multicolumn{2}{c}{$x$}
&\multicolumn{2}{c}{$y$} 
&\multicolumn{2}{c}{$E$} 
&\multicolumn{2}{c}{$\frac{1}{6}\left(\sum_{i=1}^nx_i^3\right)$}\\\midrule

$  1$&$    14$&$    -7$&$7220$&$    -0$&$1627$&$    20$&$0000$&$    89$&$9136$\\
$  1$&$    14$&$     6$&$3887$&$    -1$&$1706$&$    32$&$7773$&$    206$&$3830$\\
$  2$&$    13$&$    -5$&$4741$&$     0$&$0729$&$    21$&$8963$&$    111$&$9900$\\
$  2$&$    13$&$     4$&$1407$&$    -1$&$4063$&$    36$&$5629$&$    184$&$3060$\\
$  3$&$    12$&$    -4$&$4379$&$     0$&$2761$&$    26$&$6274$&$    123$&$0070$\\
$  3$&$    12$&$     3$&$1046$&$    -1$&$6095$&$    38$&$6274$&$    173$&$2900$\\
$  4$&$    11$&$    -3$&$7936$&$     0$&$4704$&$    30$&$3489$&$    130$&$4600$\\
$  4$&$    11$&$     2$&$4603$&$    -1$&$8037$&$    39$&$6822$&$    165$&$8360$\\
$  5$&$    10$&$    -3$&$3333$&$     0$&$6667$&$    33$&$3333$&$    136$&$2960$\\
$  5$&$    10$&$     2$&$0000$&$    -2$&$0000$&$    40$&$0000$&$    160$&$0000$\\
$  6$&$     9$&$    -2$&$9761$&$     0$&$8729$&$    35$&$7128$&$    141$&$3050$\\
$  6$&$     9$&$     1$&$6427$&$    -2$&$2063$&$    39$&$7128$&$    154$&$9910$\\
$  7$&$     8$&$    -2$&$6825$&$     1$&$0972$&$    37$&$5547$&$    145$&$9080$\\
$  7$&$     8$&$     1$&$3491$&$    -2$&$4305$&$    38$&$8880$&$
150$&$3880$\\\bottomrule
\end{tabular}
\end{center}

The value $p=5$ corresponds to
the maximal energy $E=40$ and we get the Clebsch graph with the spectrum
$10, 2^5, -2^{10}$. There are no other graphs with this spectrum
\cite{CvRS4}.

\bigskip
\noindent {\bf B.}~ Consider now graphs on $n=10$ vertices.

\medskip
1. $m=30$ edges.

\medskip

A maximal energy graph in $\cal G$ cannot be complete. Therefore it
must have at least three distinct eigenvalues.  Restrict our search
further to connected graphs with largest eigenvalue 6 (which is
simple).  Then ${\cal K}=\{6\}$.  (The complement of the Petersen
graph is in $\cal G$.)  Now we have $\setcard{K}=1$, $|J|=9$, $C=6$, and
$D=36$. Our program gives the following solutions.

\begin{center}
\begin{tabular}{rrr@{.}lr@{.}lr@{.}lr@{.}l}\toprule
 \multicolumn{1}{c}{$p$} 
&\multicolumn{1}{c}{$q$}
&\multicolumn{2}{c}{$x$}
&\multicolumn{2}{c}{$y$} 
&\multicolumn{2}{c}{$E$} 
&\multicolumn{2}{c}{$\frac{1}{6}\left(\sum_{i=1}^nx_i^3\right)$}\\\midrule
$  1$  &  $ 8$  &  $-4$&$8830$  &  $-0$&$1396$  &  $12$&$0000$  &  $16$&$5911$\\
$  1$  &  $ 8$  &  $ 3$&$5497$  &  $-1$&$1937$  &  $19$&$0994$  &  $41$&$1866$\\
$  2$  &  $ 7$  &  $-3$&$4555$  &  $ 0$&$1302$  &  $13$&$8221$  &  $22$&$2487$\\
$  2$  &  $ 7$  &  $ 2$&$1222$  &  $-1$&$4635$  &  $20$&$4888$  &  $35$&$5290$\\
$  3$  &  $ 6$  &  $-2$&$7749$  &  $ 0$&$3874$  &  $16$&$6491$  &  $25$&$3752$\\
$  3$  &  $ 6$  &  $ 1$&$4415$  &  $-1$&$7208$  &  $20$&$6491$  &  $32$&$4025$\\
$  4$  &  $ 5$  &  $-2$&$3333$  &  $ 0$&$6667$  &  $18$&$6667$  &  $27$&$7778$\\
$  4$  &  $ 5$  &  $ 1$&$0000$  &  $-2$&$0000$  &  $20$&$0000$  &  $30$&$0000$\\\bottomrule
\end{tabular}
\end{center}

High values of energy are obtained for $p=3$ and $p=2$, but the
corresponding graphs do not exist since the third spectral moment is
not divisible by 6 (cf. Proposition~1). For $p=4$, we find that
$L(K_5)$, the line graph of the complete graph $K_5$, has distinct
eigenvalues 6, 1, $-2$ with multiplicities 1, 4, 5 respectively.  The
graph $L(K_5)$ is known to have maximal energy among all 10 vertex
graphs \cite{CCGH}.

One should be aware that our calculations do not prove that $L(K_5)$
has maximal energy (this fact we know from \cite{CCGH}).  If a graph
existed that had spectrum $6$, $1.4415^3$, $-1.7208^6$, then Theorem~1
would establish it as a maximal energy graph over $\cal G$.  However,
we do not know if there is a graph with, say, eight distinct
eigenvalues (all but one close to $1.4415$ or $-1.7208$) which could
have energy very close 20.6491.  It seems that one has to introduce a
kind of distance between spectra (families of reals) to handle such
effects, but we shall not do that in this paper.

\vspace{5mm}

2. $m=9$ edges.

\medskip 

Restrict our search to graphs having a simple eigenvalue
${\cal K}=\{3\}$.  Thus $\setcard{K}=1$, $|J|=9$, $C=3$, and $D=9$.  We then
get the following solutions.
\begin{center}
\begin{tabular}{rrr@{.}lr@{.}lr@{.}lr@{.}l}\toprule
 \multicolumn{1}{c}{$p$} 
&\multicolumn{1}{c}{$q$}
&\multicolumn{2}{c}{$x$}
&\multicolumn{2}{c}{$y$} 
&\multicolumn{2}{c}{$E$} 
&\multicolumn{2}{c}{$\frac{1}{6}\left(\sum_{i=1}^nx_i^3\right)$}\\\midrule
$   1   $&$   8   $&$    -3$&$0000   $&$     0$&$0000   $&$   6$&$0000   $&$  0$&$0000   $\\
$   1   $&$   8   $&$     2$&$3333   $&$    -0$&$6667   $&$  10$&$6667   $&$  6$&$2222   $\\
$   2   $&$   7   $&$    -2$&$0972   $&$     0$&$1706   $&$   8$&$3887   $&$  1$&$4313   $\\
$   2   $&$   7   $&$     1$&$4305   $&$    -0$&$8373   $&$  11$&$7220   $&$  4$&$7910   $\\
$   3   $&$   6   $&$    -1$&$6667   $&$     0$&$3333   $&$  10$&$0000   $&$  2$&$2222   $\\
$   3   $&$   6   $&$     1$&$0000   $&$    -1$&$0000   $&$  12$&$0000   $&$  4$&$0000   $\\
$   4   $&$   5   $&$    -1$&$3874   $&$     0$&$5099   $&$  11$&$0994   $&$  2$&$8300   $\\
$   4   $&$   5   $&$     0$&$7208   $&$    -1$&$1766   $&$  11$&$7661   $&$  3$&$3922   $\\\bottomrule
\end{tabular}
\end{center}

For $p=1$ we get the star $K_{1,9}$ (in accordance with the known fact that
stars are minimal energy graphs among the trees with a fixed number of vertices)
and for $p=3$, we get that the graph
$K_4 \cup 3K_2$ has maximal energy (among graphs considered).

\bigskip
\noindent {\bf C.}~ The procedure developed is very flexible in the
sense that some structural restrictions on graphs (e.g., regularity,
bipartiteness, connectedness, etc.) can easily be introduced. Consider
graphs on $n=14$ vertices with $m=21$ edges.  If we apply Theorem 1
with ${\cal K} = \emptyset$ (general graphs), and ${\cal K} = \{3\}$
(cubic graphs, cf. Proposition~3), we do not get any
graphs. However, if we assume ${\cal K} = \{3, -3\}$ (cubic
bipartite graphs, cf. Proposition~2), we obtain the Heawood graph with
eigenvalues $3, \sqrt{2}^6, -\sqrt{2}^6, -3$, which really is an
extremal graph.

\section{Some data on maximal energy graphs}

The presented material might suggest that graphs with extremal
energy always have a small number of distinct eigenvalues. We
now present some computational and theoretical facts showing
that this is not always the case.

We present in Table~\ref{tab:max} results of computations concerning extremal energy
graphs up to $n=12$ vertices, performed when the paper \cite{CCGH} was
being prepared, but which have not been published.   Each graph is
identified with a graph6 code, which is a compact representation of
the adjacency matrix.  The specification of the graph6 code is
distributed with Brendan McKay's Nauty program \cite{McK} and can also
be found on the Nauty website.  We have independently verified these
results with our own computer programs by an exhaustive search.
For $n=11, 12$ we reduced the search to specific numbers of
edges ($m=35, 36, 37$ for $n=11$ and $m=41, 42, 43$ for $n=12$---about
4.6 billion graphs in the latter case).

\begin{table}
  \centering
  \begin{tabular}{cccccp{4cm}}\toprule
    Vertices & Graph6  & Edges & Energy & \parbox{2cm}{Distinct \\Eigenvalues} & Spectrum\\\midrule
    7 &\verb@F`~~w@& 17 &12 & 4 &$5, 1, -1^4, -2$\\
    8 &\verb@G`lv~{@& 21 & $14.325$ & 7 & $5.427$, $1.118$, $0.618$, $-1^2$, $-1.618$, $-1.679$, $-1.865$\\
    9 &\verb@HEutZhj@ & 21 & 17.060 & 6 & $4.702$, $1.414^2$, $1$, $-1.414^2$, $-1.702$, $-2^2$\\
    10 & \verb@I~qkzXZLw@ & 30 & 20 & 3 & $6, 1^4, -2^5$\\
    11  &\verb@JJ^em]uj[v_@& 36 & 22.918 & 5 & $6.585$, $1.874$, $1^3$,
$-1.459$, $-2^5$\\
12 & \verb@K~z\c\qRXVa~@ & 42 & 26 & 5 & $7$, $2^2$, $1^2$, $-1$, $-2^6$\\\bottomrule
  \end{tabular}
  \caption{Maximal energy graphs}
  \label{tab:max}
\end{table}

We note a few items of interest from Table~\ref{tab:max}.  For
$n=7$ there exists a graph with the spectrum $5, 1, -1^4, -2$ which
has the same (maximal) energy ($E=12$) as the complete graph $K_7$.
For $n=12$ the extremal graph is regular of degree 7.  This graph is
known as an exceptional graph for the least eigenvalue $-2$ and can be
found as graph No.~186 in Table~A3 of~\cite{CvRS4}.  In these
examples, the eigenvalue $-2$ almost always appears with high
multiplicity. The number $-2$ is the least eigenvalue in almost all
line graphs. For the role of the number $-2$ in the theory of graph
spectra, see~\cite{CvRS4}.

This data partially confirms the tendency of eigenvalues of high
multiplicity (i.e., small number of distinct eigenvalues) in maximal
energy graphs.  However, some cases in which there are a large number
of distinct eigenvalues also occur.  Another example of graphs with
maximal energy and a high number of distinct eigenvalues is known from
theoretical considerations.  Among trees with a given number of
vertices, the path has maximal energy (cf., e.g., \cite{CvDSa},
p.~238).

Examples with a high number of distinct eigenvalues show that, in
these cases, several sets of eigenvalues produced by procedure $L$ by
specifying families of eigenvalues $\cal K$ do not actually correspond
to graphs.  This can be seen in the case of a path in the following
way.  Let $\pi_1,\pi_2,\dots,\pi_n$ be (distinct) eigenvalues of a path
$P_n$ on $n$ vertices. Let ${\cal K}_i = \{\pi_1,\pi_2,\dots,\pi_i\},
i=1,2,\dots,n-2$.  By applying procedure $L$ in turn with these sets we
come to $P_n$ only after $n-2$ iterations.

\section{A heuristic procedure}

Having in view all that has been said, one can outline the
following heuristic procedure for finding extremal graphs by
repeated applications of procedure~$L$.

\medskip\noindent
{\bf Procedure.} {\em Let $n,m$ be positive integers and let
  $\cal K$ be the family of reals.  Let
  ${\cal K}_L$ be the eigenvalues derived from procedure $L$.  Try to
  construct a graph having eigenvalues ${\cal K}_L$.  If such a graph
  does not exist, repeatedly add eigenvalues to $\cal K$ and try to
  construct a graph with eigenvalues ${\cal K}_L$.}

\medskip

Of course, depending on the concrete problem, the eigenvalues to add
to $\cal K$ should be determined with other facts and tools not
contained in Theorem~1.

We shall give an example of using this procedure to determine a
regular graph of degree 15 on 18 vertices with maximal energy.  Since
the complement of such a graph is regular of degree 2, we know that
the complement of such a graph is composed of disjoint unions of
cycles.    Hence we could apply our Mathematica program a
few times while trying, for example, the presence of eigenvalues $-\varphi \approx -1.618,
1/\varphi \approx 0.618$, which come from a pentagon in the complement, and/or $1, -1$,
which come from a quadrangle in the complement (cf. Propositions~4 and
5).  However, we will try to use procedure $L$ to guide our guesses.

Using our Mathematica program with $n=18$, $m=15n/2=135$ and ${\cal
  K}=\{15\}$, we have the following solutions:
\begin{center}
\begin{tabular}{rrr@{.}lr@{.}lr@{.}lr@{.}l}\toprule
 \multicolumn{1}{c}{$p$} 
&\multicolumn{1}{c}{$q$}
&\multicolumn{2}{c}{$x$}
&\multicolumn{2}{c}{$y$} 
&\multicolumn{2}{c}{$E$} 
&\multicolumn{2}{c}{$\frac{1}{6}\left(\sum_{i=1}^nx_i^3\right)$}\\\midrule
$   1$&    $16$&    $-6$&$3501$&    $-0$&$5406$&    $30$&$0000$&    $519$&$4020$\\
$   1$&    $16$&    $ 4$&$5854$&    $-1$&$2241$&    $39$&$1708$&    $573$&$6770$\\
$   2$&    $15$&    $-4$&$6259$&    $-0$&$3832$&    $30$&$0000$&    $529$&$3640$\\
$   2$&    $15$&    $ 2$&$8612$&    $-1$&$3815$&    $41$&$4446$&    $563$&$7160$\\
$   3$&    $14$&    $-3$&$8353$&    $-0$&$2496$&    $30$&$0000$&    $534$&$2570$\\
$   3$&    $14$&    $ 2$&$0706$&    $-1$&$5151$&    $42$&$4234$&    $558$&$8230$\\
$   4$&    $13$&    $-3$&$3466$&    $-0$&$1241$&    $30$&$0000$&    $537$&$5080$\\
$   4$&    $13$&    $ 1$&$5819$&    $-1$&$6406$&    $42$&$6554$&    $555$&$5720$\\
$   5$&    $12$&    $-3$&$0000$&    $ 0$&$0000$&    $30$&$0000$&    $540$&$0000$\\
$   5$&    $12$&    $ 1$&$2353$&    $-1$&$7647$&    $42$&$3529$&    $553$&$0800$\\
$   6$&    $11$&    $-2$&$7332$&    $ 0$&$1272$&    $32$&$7983$&    $542$&$0860$\\
$   6$&    $11$&    $ 0$&$9685$&    $-1$&$8919$&    $41$&$6218$&    $550$&$9940$\\
$   7$&    $10$&    $-2$&$5162$&    $ 0$&$2613$&    $35$&$2261$&    $543$&$9450$\\
$   7$&    $10$&    $ 0$&$7514$&    $-2$&$0260$&    $40$&$5203$&    $549$&$1350$\\
$   8$&    $ 9$&    $-2$&$3322$&    $ 0$&$4064$&    $37$&$3153$&    $545$&$6870$\\
$   8$&    $ 9$&    $ 0$&$5675$&    $-2$&$1711$&    $39$&$0800$&    $547$&$3930$\\\bottomrule
\end{tabular}
\end{center}

As we can see, the only solution that passes the 3rd moment test is
$p=5$, $q=12$, $x=-3$, $y=0$, and $E=30$.  Interestingly, this is the
lowest energy solution found. We have obtained the eigenvalues of
the graph $\overline{6C_3}$ with the spectrum $15, 0^{12},-3^5$. It
has three distinct eigenvalues and belongs to the set of trivial
strongly regular graphs.

By Propositions~3 and 4, a 15-regular 18-vertex graph has eigenvalue
$-3$ with multiplicity $k$ if and only if the complement of the graph
has $k+1$ components (which are all cycles).  Let us assume for the
moment that our extremal graph has at least 4 components in its
complement. Therefore we take ${\cal K}=\{15,-3,-3,-3\}$. We have the
following computer output (in the tables below, the lines in which the third
spectral moment is divisible by 6 will be denoted by $+$ and the others
will be denoted by $-$).

\begin{center}
\begin{tabular}{rrr@{.}lr@{.}lr@{.}lc}\toprule
 \multicolumn{1}{c}{$p$} 
&\multicolumn{1}{c}{$q$}
&\multicolumn{2}{c}{$x$}
&\multicolumn{2}{c}{$y$} 
&\multicolumn{2}{c}{$E$} 
&\multicolumn{1}{c}{3rd moment test}\\\midrule
$1$&$13$& $-4$&$2136$&$ -0$&$1374$&$ 30$&$0000$&$- $\\
$1$&$13$& $ 3$&$3565$&$ -0$&$7197$&$ 36$&$7129$&$- $\\
$2$&$12$& $-3$&$0000$&$  0$&$0000$&$ 30$&$0000$&$+ $\\
$2$&$12$& $ 2$&$1429$&$ -0$&$8571$&$ 38$&$5714$&$- $\\
$3$&$11$& $-2$&$4387$&$  0$&$1197$&$ 32$&$6325$&$- $\\
$3$&$11$& $ 1$&$5816$&$ -0$&$9768$&$ 39$&$4896$&$- $\\
$4$&$10$& $-2$&$0884$&$  0$&$2354$&$ 34$&$7074$&$- $\\
$4$&$10$& $ 1$&$2313$&$ -1$&$0925$&$ 39$&$8502$&$- $\\
$5$&$ 9$& $-1$&$8370$&$  0$&$3539$&$ 36$&$3700$&$- $\\
$5$&$ 9$& $ 0$&$9799$&$ -1$&$2110$&$ 39$&$7986$&$- $\\
$6$&$ 8$& $-1$&$6408$&$  0$&$4806$&$ 37$&$6891$&$- $\\
$6$&$ 8$& $ 0$&$7836$&$ -1$&$3377$&$ 39$&$4033$&$- $\\
$7$&$ 7$& $-1$&$4784$&$  0$&$6212$&$ 38$&$6969$&$- $\\
$7$&$ 7$& $ 0$&$6212$&$ -1$&$4784$&$ 38$&$6969$&$-$\\\bottomrule
\end{tabular}
\end{center}

Note that third line passes the 3rd moment test and yields again
$\overline{6C_3}$.

The second line for $p=6$ shows that we can get high energy
solutions if we have a group of six eigenvalues around  0.7836 and
a group of eight eigenvalues around  $-1.3377$. This can be
roughly achieved if we assume the existence of two pentagons in
the complement.  By Proposition~4, each pentagon in the complement
will introduce two eigenvalues of $\varphi-1\approx 0.618$ and two eigenvalues of
$-1/\varphi-1=-\varphi\approx -1.618$.  Therefore we take
$${\cal K}=\{15, -3, -3, -3,  \varphi-1, \varphi-1, \varphi-1, \varphi-1,
-\varphi, -\varphi, -\varphi, -\varphi\}$$
which yields

\begin{center}
\begin{tabular}{rrr@{.}lr@{.}lr@{.}lc}\toprule
 \multicolumn{1}{c}{$p$} 
&\multicolumn{1}{c}{$q$}
&\multicolumn{2}{c}{$x$}
&\multicolumn{2}{c}{$y$} 
&\multicolumn{2}{c}{$E$} 
&\multicolumn{1}{c}{3rd moment test}\\\midrule
$ 1$ & $ 5$ & $-2$&$4415$ & $ 0$&$0883$ & $35$&$8273$ & $-$\\
$ 1$ & $ 5$ & $ 1$&$7749$ & $-0$&$7550$ & $38$&$4940$ & $-$\\
$ 2$ & $ 4$ & $-1$&$6667$ & $ 0$&$3333$ & $37$&$6109$ & $-$\\
$ 2$ & $ 4$ & $ 1$&$0000$ & $-1$&$0000$ & $38$&$9443$ & $+$\\
$ 3$ & $ 3$ & $-1$&$2761$ & $ 0$&$6095$ & $38$&$6011$ & $-$\\
$ 3$ & $ 3$ & $ 0$&$6095$ & $-1$&$2761$ & $38$&$6011$ & $-$\\\bottomrule
\end{tabular}
\end{center}

Note that the fourth line:
\begin{center}
\begin{tabular}{rrr@{.}lr@{.}lr@{.}lc}
$ 2$ & $ 4$ & $ 1$&$0000$ & $-1$&$0000$ & $38$&$9443$ & $+$\\
\end{tabular}
\end{center}
passes the 3rd moment test.  This spectrum is realized by the
complement of $2C_4\cup2C_5$.  Through an exhaustive search, we have
verified that this graph has maximal energy among all regular graphs
of degree 15 on 18 vertices.  This extremal graph,
$\overline{2C_4\cup2C_5}$, has 6 distinct eigenvalues and spectrum $15,
1^2, (\varphi-1)^4, -1^4, -\varphi^4, -3^3$ (cf. Propositions~4 and 5).

We could have approached this another way. If we assume the
presence of two quadrangles by letting
$${\cal K}=\{15,-3,-3,-3,  -1, -1, -1, -1, 1, 1\},$$
we shall get the same solution, as the
following output shows.
\begin{center}
\begin{tabular}{rrr@{.}lr@{.}lr@{.}lc}\toprule
 \multicolumn{1}{c}{$p$} 
&\multicolumn{1}{c}{$q$}
&\multicolumn{2}{c}{$x$}
&\multicolumn{2}{c}{$y$} 
&\multicolumn{2}{c}{$E$} 
&\multicolumn{1}{c}{3rd moment test}\\\midrule
 $1$ & $  7$ & $ -3$&$4580$ & $ -0$&$0774$ & $ 34$&$0000$ & $ -$\\
 $1$ & $  7$ & $  2$&$4580$ & $ -0$&$9226$ & $ 38$&$9161$ & $ -$\\
 $2$ & $  6$ & $ -2$&$4365$ & $  0$&$1455$ & $ 35$&$7460$ & $ -$\\
 $2$ & $  6$ & $  1$&$4365$ & $ -1$&$1455$ & $ 39$&$7460$ & $ -$\\
 $3$ & $  5$ & $ -1$&$9434$ & $  0$&$3660$ & $ 37$&$6603$ & $ -$\\
 $3$ & $  5$ & $  0$&$9434$ & $ -1$&$3660$ & $ 39$&$6603$ & $ -$\\
 $4$ & $  4$ & $ -1$&$6180$ & $  0$&$6180$ & $ 38$&$9443$ & $ +$\\
 $4$ & $  4$ & $  0$&$6180$ & $ -1$&$6180$ & $ 38$&$9443$ & $ +$\\\bottomrule
\end{tabular}
\end{center}
Note that the last two (equivalent) lines pass the 3rd moment
test. Although we get the same solution, the conclusion is less
clear since we have solutions with higher energy.

\section{Conclusion}

The material presented in this paper shows how one might look for
extremal graphs using tools from calculus. According to this approach,
graphs with extremal energy ``should'' have a small number of
distinct eigenvalues.  However, the discrete
nature of the problem often prevents the expected ``nice'' solutions
from existing.  In many cases, candidates other than the theoretically
optimal solution must be considered.  However, even if the techniques
in this paper do not directly yield an extremal graph, they can also
be used to guide a search for an extremal graph.

One can say that results by J.H.~Koolen and V.~Moulton \cite{KM} and
V.~Nikiforov \cite{Nik} give a solution of the maximal energy
problem in an asymptotic sense. As far as we know at the present,
a small number of distinct eigenvalues in extremal graphs appears rarely
(complete graphs for $n \leq 7$ and strongly regular graphs for $n=10$
and $n=4\tau^2$). Our results for $n=8,9,10,11,12$ show
that, at least for moderate values of $n$, the structure of maximal
graphs could vary unexpectedly. Having in mind the fact that paths are
maximal energy trees, such effects also appear in some form for large
values of $n$.

It remains open whether it will be possible to describe maximal
energy graphs in a better way using new approaches.

\renewcommand{\refname}{Bibliography}

{\footnotesize \flushleft
D.~Cvetkovi\' c\\
Faculty of Electrical Engineering,\\
University of Belgrade,\\
P.O.Box 35--54,\\
11120 Belgrade, Serbia\\
e-mail: ecvetkod@etf.bg.ac.yu\\
 }

\medskip
{\footnotesize \flushleft
J.~Grout \\
Department of Mathematics \\
Brigham Young University \\
Provo, UT 84602 \\
USA \\
e-mail: grout@math.byu.edu\\ }

\bigskip

\appendix
\section{Mathematica Program}
\label{app:mathematica-program}
\begin{verbatim}
possibleEvals[numvertices_Integer, numedges_Integer, knownevals_List] :=
 Block[{eq, p, q, x, y, en},
   eq := {p + q == numvertices - Length[knownevals],
       p*x + q*y == -(Plus @@ knownevals),
       p*x^2 + q*y^2 == 
      2*numedges - (Plus @@ (knownevals^2))}; Print[ InputForm[eq]];
   Table[{p, q, x, y,  p*Abs[x] + q*Abs[y] + Plus @@ (
              Abs[knownevals]), (p*x^3 + q*
                y^3 + (Plus @@ (knownevals^3)))/6} /.
         Solve[eq // Append[#, p == i] &, {p, q, x, y}], {i, 1, 
            Ceiling[(numvertices - Length[knownevals] - 1)/2]}] //
      Flatten[#, 1] &]
\end{verbatim}
\section{SAGE Program}
\label{app:sage-program}
\begin{verbatim}
def possible_evals(num_vertices, num_edges, known_evals):
    # Make x and y variables
    var('x, y')
    upper_bound=ceil((1/2)*(num_vertices - len(known_evals) - 1))
    for p in [1..upper_bound]:
        q = num_vertices - len(known_evals) - p
        solutions = solve([p*x + q*y == -sum(known_evals), \
          p*x^2 + q*y^2 == 2*num_edges - sum([i^2 for i in known_evals])], \
          x,y, solution_dict=True)
        for s in solutions:
            energy = p*abs(s[x]) + q*abs(s[y]) \
                     + sum([abs(e) for e in known_evals])
            third_moment_test = (1/6)*(p*s[x]^3 + q*s[y]^3 \
                                       + sum([e^3 for e in known_evals]))
            print "%d, %d, %f, %f, %f, %f"%(p, q, s[x], s[y], \
                                            energy, third_moment_test)
\end{verbatim}

\end{document}